\documentclass[a4paper]{article}

\usepackage{xcolor}

\newcommand{\vsig}{\mbox{\boldmath $\sigma$}}

\newcommand{\vx}{\mbox{\boldmath $x$}}
\newcommand{\vdel}{\mbox{\boldmath $\delta$}}

\newcommand{\had}{\bullet}

\newcommand{\ft}{\hat}

\newcommand{\msgguy}[1]{{\textcolor{black}{#1}}}
\newcommand{\corguy}[1]{{\textcolor{black}{#1}}}

\newcommand{\hauteur}{5cm}

\input cyracc.def
\font\tencyr=wncysc10
\def\cyr{\tencyr\cyracc}
\def\diracComb{\mbox{\cyr SH}}

 \DeclareFontFamily{U}{wncy}{}
    \DeclareFontShape{U}{wncy}{m}{n}{<->wncyr10}{}
    \DeclareSymbolFont{mcy}{U}{wncy}{m}{n}
    \DeclareMathSymbol{\Sh}{\mathord}{mcy}{"58} 

\usepackage{amsmath,amssymb}
\usepackage{graphicx,epstopdf}
\usepackage[caption=false]{subfig}

\DeclareMathOperator{\sinc}{sinc}

\usepackage[pages=all, color=black, position={current page.south}, placement=bottom, scale=1, opacity=1, vshift=5mm]{background}
\SetBgContents{
	\tt This work is shared under a \href{https://creativecommons.org/licenses/by-sa/4.0/}{CC BY-NC-ND 4.0 license} unless otherwise noted
}      

\usepackage[margin=1in]{geometry} 

\usepackage{amsthm}

\usepackage[utf8]{inputenc}
\usepackage{hyperref}
\hypersetup{
	unicode,
	pdfauthor={Author One, Author Two, Author Three},
	pdftitle={A simple article template},
	pdfsubject={A simple article template},
	pdfkeywords={article, template, simple},
	pdfproducer={LaTeX},
	pdfcreator={pdflatex}
}

\usepackage[sort&compress,numbers,square]{natbib}
\theoremstyle{plain}

\theoremstyle{definition}

\usepackage{graphicx, color}
\graphicspath{{fig/}}

\usepackage{algorithm, algpseudocode} 
\usepackage{mathrsfs} 

\usepackage{lipsum}

\title{An exact solution to the Fourier Transform of band-limited periodic functions with  nonequispaced data and application to non-periodic functions}
\author{Guy Perrin$^1$}

\date{
	$^1$LESIA, Observatoire de Paris-PSL, CNRS, Sorbonne Universit\'e, Universit\'e de Paris, 92190 Meudon, France \\ \texttt{guy.perrin@obspm.fr}\\%
}

\begin{document}
\maketitle





\begin{abstract}
The need to Fourier transform data sets with irregular sampling is shared by various domains of science. This is the case for example in astronomy or sismology. Iterative methods have been developed that allow to reach approximate solutions. Here an exact solution to the problem for band-limited periodic signals is presented. The exact spectrum can be deduced from the spectrum of the non-equispaced data through the inversion of a Toeplitz matrix. The result applies to data of any dimension. This method also provides an excellent approximation for non-periodic band-limit signals. The method allows to reach very high dynamic ranges ($10^{13}$ with double-float precision) which depend on the regularity of the samples. 
\end{abstract}




 
\section{Introduction}
\label{section:introduction}
The issue of nonuniform sampling is quite common in various fields of science or engineering. It arises for example in the study of variable astronomical sources because they are not observable all year round or their observability may be impaired by adverse weather conditions. It is also a common problem in aperture synthesis with non-uniformly filled spatial frequency planes. Other examples can be found in sismology, tomography, \dots. In some cases, the scientific measurement is derived from the Fourier spectrum. Several issues are the consequence of the nonuniform sampling of the signal. First, the efficient algorithms to compute discrete Fourier spectra as the Fast Fourier Transform cannot be used. Second, spurious features are produced in the reconstructed spectrum that lead to a reduction of signal-to-noise ratio and plague the detection of faint signals. Solutions have been proposed to efficiently and directly compute the spectrum of non-equispaced data~\cite{Dutt1993}. The reconstruction of spectra is equivalent to the reconstruction of the signal as one is deduced from the other by a Fourier transform. Iterative algorithms are known to converge to the exact signal in the band-limited case with excellent accuracy \cite{feichtinger2021}. Other methods are based on linear combinations of interpolating functions without the need for iterations \cite{Cheng2020}. These methods can be exact in the case of regular sampling or for periodic functions or for recurrent irregular sampling following the Generalized Sampling \msgguy{Expansion} scheme for the latter \cite{Papoulis1977}. But extra informations (this is the case for the Multi-Channel Interpolation method) or a high sampling rate are required to reach high accuracies in the general case. 

Here we propose a direct method to recover the spectrum of a band-limited function. The method is exact in the case of a periodic and band-limited signal and is an approximation if the function is not periodic. It was first skectched in 1996 \cite{Perrin1996} and is developed in more details here. 

The paper is organized as follows. The formalism of irregular sampling and the link between the spectrum of the nonequispaced data and the Fourier spectrum are established in Section~\ref{sec:sampling}. The exact solution of the Fourier transform of nonuniform samples is presented in Section~\ref{sec:spectrum_periodic} for periodic band-limited functions while the approximated solution for non-periodic but band-limited functions is discussed in Section~\ref{sec:spectrum_non-periodic}. The complexity and accuracy of the algorithm is detailed in Section~\ref{sec:accuracy}. \msgguy{The method is compared with other algorithms in Section~\ref{section:comparison}}. Last, the extension of the method to dimensions higher than 1 is derived in Section~\ref{sec:dimension} before conclusions are given in Section~\ref{sec:conclusion}. 
 
\section{Irregular sampling}
\label{sec:sampling}
\subsection{Definitions}

$S$ is a band-limited signal and depends on the variable $x$ defined over $\Bbb{R}$. $S$ takes values in $\Bbb{C}$ and it is assumed that $S \in {L}^2(\Bbb{R})$. Its spectrum is noted $\ft{S}$ and is the Fourier transform of $S$ where we choose the form of the Fourier transform where the variable conjugated to $x$ is the associated frequency $\sigma$. If $x$ is time then $\sigma$ is the temporal frequency expressed in Hertz. If $x$ is the optical path difference then $\sigma$ is the wavenumber classically expressed in cm$^{-1}$ by spectroscopists. $\ft{S}$ writes:
\begin{equation}
\ft{S}(\sigma)=\int_{-\infty}^{+\infty}{S(x)\, e^{-2i\pi\sigma x}\, dx}
\end{equation}

\noindent $\ft{S}$ takes values in $\Bbb{C}$ and is zero outside the interval $[-f_\mathrm{max},+f_\mathrm{max}]$. We assume $S$ is sampled over the grid $\Gamma=\{ x_{0},\ldots,x_{N-1} \}$ with $x_{0} < \ldots < x_{N-1}$. The sampled version of $S$ is noted $S_{\Gamma}$. It is the sum of Dirac distributions centered on the points of the grid. In this paper, distributions are noted as functions rather than as linear forms over the space of functions. Following this principle, $S_{\Gamma}$ writes: 
\begin{equation} 
S_{\Gamma}(x)=\sum_{k=0}^{N-1}{S(x_{k})\,\delta (x-x_{k}}) 
\end{equation}

\noindent The Fourier transform of $S_{\Gamma}$ writes:
\begin{equation}
\label{eq:spectrum}
\ft{S}_{\Gamma}(\sigma)=\sum_{k=0}^{N-1}{S(x_{k})\,e^{-2i\pi\sigma x_{k}}}
\end{equation}

\noindent Let us define the average sampling interval $\overline{\delta x}$:
\begin{equation}
\overline{\delta x}=\frac{x_{N-1}-x_{0}}{N-1}
\end{equation}

\noindent We associate the width $\Delta x$ to the sampling grid $\Gamma$:
\begin{equation}
\Delta x = N\, \overline{\delta x}
\end{equation}

\noindent We define a regularly sampled version $S_{\mathrm{Reg}}$ of $S$ using the grid of equispaced samples $\mathrm{Reg}=\{ x_{0}, x_{0} + \overline{\delta x}, \ldots,x_{0} +(N-1)\overline{\delta x}\}$. $S_{\Gamma}$ and $S_{\mathrm{Reg}}$ share the same first and last samples, the same sampling interval and the same width. \\ \\
{\it{Remark 1: case of periodic signals.}} In case the underlying function of signal $S$ is $P-$periodic it is necessary to choose $\Delta x$ as a multiple of $P$ so that the exact inversion of Section~\ref{sec:spectrum_periodic} holds. We then simply choose $\Delta x = P \times (x_{N-1}-x_0)//P$ and define $\overline{\delta x} = \frac{\Delta x}{N}$. \\ \\
{\it{Remark 2: case of decimated samples.}} In case the sampled distribution is a decimated version of equispaced samples, the regular samples are the undecimated samples limited to $N$ samples. In this case $\overline{\delta x} = \delta x$. 
\noindent As for regular sampling, $\Delta x$ must be strictly larger than $(x_{N-1}-x_{0})$ to allow the periodic replication of the samples to sample the spectrum (see equation~(\ref{eq:periodized_samples})). Let us define the door function $ \Pi$ of width $\Delta x$  as the indicator function of the interval $[x_{0}-\frac{\overline{\delta x}}{2},x_{0}+\Delta x
-\frac{\overline{\delta x}}{2}]$:
%
%
\begin{equation} \Pi(x)= \left\{ \begin{array}{ll}
                          1 & \mbox{if $x \in [x_{0}-\frac{\overline{\delta x}}{2},x_{0}+\Delta x
-\frac{\overline{\delta x}}{2}]$} \\
                          0 & \mbox{otherwise}
\end{array} \right. \end{equation}
%
%
%
Last, $\diracComb_{P}$ is the Dirac comb of step $P$, with $P>0$ :
\begin{equation}
\diracComb_P(x)=\sum_{n \in \mathbb{Z}}{\delta (x-nP)}
\end{equation}
whose Fourier Transform writes:
\begin{equation}
\ft\diracComb_P(\sigma)=\frac{1}{P}\sum_{n \in \mathbb{Z}}{\delta (\sigma-\frac{n}{P})}=\frac{1}{P}\diracComb_\frac{1}{P}(\sigma)
\end{equation}


\subsection{The spectrum of the sampled signal}

Let us define the sampling distribution $G_\Gamma$ derived from the grid $\Gamma$:
%
\begin{equation}
G_\Gamma(x)=\sum_{n \in \mathbb{Z}}{\delta (x-\tilde{x}_{n})}
\end{equation}
\noindent where:
\begin{equation}
\tilde{x}_{n}= \left\{ \begin{array}{ll}
                          x_{n} & \mbox{for $n \in [0,..,N-1]$} \\
                          x_{r_{n}}+(n-{r_{n}})\, {\overline{\delta x}} & \mbox{otherwise}									\end{array} \right.
\end{equation}
with  $r_{n}$ the remainder of the division of $n$ by $N$. We call $\hat{S}_\Gamma$ the signal $S$ sampled over the extended grid $\left\{ \tilde{x}_n \right\}_{n \in \Bbb{Z}}$. The extended grid is the combination of $N$ interleaved regular grids each of step $\Delta x$ which constitute a scheme of recurrent irregular sampling. We know from the generalized sampling theorem of Papoulis \cite{Papoulis1977} that $S$ can be reconstructed from the samples $\hat{S}_\Gamma$ if $\frac{\Delta x}{N} \leqslant \frac{1}{2 f_{\mathrm{max}}}$ which is equivalent to ${\overline{\delta x}}$ has to check the Nyquist sampling criterion as in the classical case of regular sampling: ${\overline{\delta x}} \leqslant \frac{1}{2 f_{\mathrm{max}}}$. \\

\noindent The extended grid infinite sum can be reorganized into a finite sum of interleaved infinite grids of step $\Delta x$:
\begin{equation}
G_\Gamma(x)=\sum_{k=0}^{N-1}{\left[ \sum_{n \in \mathbb{Z}}{\delta (x-x_{k}-n \Delta x)}
\right]}
\end{equation}
\noindent This expression can be simplified using the Dirac comb:
%
\begin{equation}
G_\Gamma(x)=\sum_{k=0}^{N-1}{ \diracComb_{\Delta x} ( x-x_{k})}
\end{equation}

\noindent The sampling grid $G_\Gamma$ is therefore a finite sum of Dirac combs whose Fourier transforms are defined. It is the underlying principle of the method: connect the irregular sampling grid to a set of periodic sampling grids. The Fourier transform of $G_\Gamma$ writes:
%
\begin{equation}
{\ft G}_\Gamma(\sigma)=\frac{1}{\Delta x} \left[ \sum_{k=0}^{N-1}{e^{-2i\pi \sigma x_{k}}}\right] 
                     \diracComb_{\delta \sigma}(\sigma)
\end{equation}

\noindent with:
\begin{equation}
\label{Eq:delta_sigma}
\Delta x \delta\sigma =1
\end{equation}

\noindent We now define the function $C_\Gamma$ as:
\begin{equation}
C_\Gamma(\sigma)=\frac{1}{\Delta x}  \sum_{k=0}^{N-1}{e^{-2i\pi \sigma x_{k}}}
\end{equation}

\noindent The Fourier transform of the sampling function $G$ can be rewritten: 
\begin{equation}
{\ft
G}_\Gamma(\sigma)=\sum_{n \in \mathbb{Z}}{C_\Gamma(n\delta\sigma)\delta(\sigma-n\delta\sigma)}
\end{equation}

\noindent The sampled signal $S_\Gamma$ can be written using the sampling function $G_\Gamma$:
\begin{equation}
\label{Eq:sampled_signal}
S_{\Gamma}(x)=S(x)\Pi (x) G_\Gamma(x)
\end{equation}

\noindent We note $\mathscr{S}=S\Pi$, the signal $S$ multiplied by the door function $\Pi$, the truncated version of $S$ over the width $\Delta x$. Though $S$ is a band-limited function, $\mathscr{S}$ is not band-limited anymore in the general case. Its sprectrum is the convolution of $\ft{S}$ by a sinc function $\frac{\sin (\pi\sigma\Delta x)}{\pi \sigma\Delta x}$ of width $\frac{1}{\Delta x}$. The spectrum of the signal sampled with the grid $\Gamma$ is the convolution of $\ft{\mathscr{S}}$ by the Fourier transform of the sampling function $G_\Gamma$:
\begin{equation}
\label{Eq:spectrum_samples}
\ft{S}_{\Gamma}(\sigma)=\ft{\mathscr{S}}(\sigma) \star {\ft G}_\Gamma(\sigma) =
\sum_{n \in \mathbb{Z}}{C_\Gamma(\sigma_{n})\ft{\mathscr{S}}(\sigma-\sigma_{n})}
\end{equation}
where $\sigma_{n}=n\delta\sigma$ for $n$  $\in \Bbb{Z}$. Defining $N_0=-\frac {N}{2}$ if $N$ is even and $-\frac {N-1}{2}$ if $N$ is odd, we call $\ft{{\cal X}} = [ \ft{X}(\sigma_{N_0}), \dots, \ft{X}(\sigma_{N_0+N-1}) ]$ the vector associated to the spectrum $\ft{X}$.
%
%
The effect of irregular sampling is to introduce a linear relation between the frequency components of the spectrum. The goal of the paper is to recover the spectrum of $S$ sampled over a regular grid, $\ft{\cal S}_{\mathrm{Reg}}$, from $\ft{\cal S}_{\Gamma}$. The periodic and non-periodic cases are respectively addressed in Section~\ref{sec:spectrum_periodic} and in Section~\ref{sec:spectrum_non-periodic}. \\ \\
{\it{Lemma: relation between the spectrum of the irregular samples and the spectrum of the regularly sampled signal.}} There exists a matrix $A$ such that: $\ft{\cal S}_{\Gamma} = A.\ft{\cal S}_{\mathrm{Reg}} +{\cal R}$, where the elements of ${\cal R}$ depend on samples of the spectrum of the truncated signal $\ft{\mathscr{S}}(\sigma_n)$ for $n$ outside the interval $[N_0,\dots,N_0+N-1]$.  \\ \\
{\it{Proof}}\\
The sampling of the spectrum of equation~(\ref{Eq:spectrum_samples}) with the sampling interval $\delta \sigma$ defined in equation~(\ref{Eq:delta_sigma}) is equivalent to replicating the sampled function of equation~(\ref{Eq:sampled_signal}) with the period $\Delta x$ leading to the $\Delta x$-periodic distribution:
%
\begin{equation}
\label{eq:periodized_samples}
S_{\Gamma}(x)  \star \Delta x \diracComb_{\Delta x}(x)
\end{equation}

\noindent The spectrum of the above signal now equates the sampled version of the spectrum of the irregularly sampled $S$ of equation~(\ref{Eq:spectrum_samples}) with sampling interval $\delta \sigma$:
\begin{equation}
\ft{S}_{\Gamma}(\sigma_{n})=\sum_{k \in \mathbb{Z}}{C_\Gamma(\sigma_{n-k})\ft{\mathscr{S}}(\sigma_{k})}
\label{eq:sampled_spectrum}
\end{equation}
It is $\Delta \sigma$-periodic with $\Delta \sigma = N \delta \sigma$.
%
%
We rewrite equation~(\ref{eq:sampled_spectrum}) by grouping terms every multiples of $N$:
\begin{equation}
\begin{split}
\ft{S}_{\Gamma}(\sigma_{n}) & =\sum_{k \in \mathbb{Z}}{C_\Gamma(\sigma_{n-k})\ft{\mathscr{S}}(\sigma_{k})} \\
& =\sum_{k \in \mathbb{Z}}
\left(\sum_{p=N_0}^{N_0+N-1}
{C_\Gamma(\sigma_{n-kN-p})\ft{\mathscr{S}}(\sigma_{kN+p})}
\right)
\end{split}
\label{eq:sampled_spectrum_general}
\end{equation}
%
The samples $\ft{S}_{\mathrm{Reg}}(\sigma_n)$ of the spectrum of the regularly sampled signal $S$ write:
\begin{equation}
\begin{split}
\ft{S}_{\mathrm{Reg}}(\sigma_{n}) & =N{\Delta x}^{-1}\sum_{k \in \mathbb{Z}}{\ft{\mathscr{S}}(\sigma_{kN+n})}\\
& = N{\Delta x}^{-1}{\ft{\mathscr{S}}(\sigma_{n})} + N{\Delta x}^{-1}\sum_{k \in \mathbb{Z}^{\star}}{\ft{\mathscr{S}}(\sigma_{kN+n})}\\
& = N{\Delta x}^{-1}{\ft{\mathscr{S}}(\sigma_{n})} + E_{\mathrm{Reg}}(\sigma_{n})
\end{split}
\label{eq:sampled_spectrum_regular}
\end{equation}
In this case, $C_{\mathrm{Reg}}(\sigma_n)=N{\Delta x}^{-1}$ if $n$ is a multiple of $N$ and 0 otherwise.
The definition of $E_{\mathrm{Reg}}(\sigma_{n})$ in the above expression will become obvious when discussing the case of periodic signals sampled over a multiple of the period. At this stage it accounts for the leakage of frequencies other than $\sigma_n$.
In the general case, $\ft{S}_{\Gamma}$ writes likewise:
\begin{equation}
\begin{split}
\ft{S}_{\Gamma}(\sigma_{n}) & = \sum_{p=N_0}^{N_0+N-1}
{C_\Gamma(\sigma_{n-p})\ft{\mathscr{S}}(\sigma_{p})} + \sum_{p=N_0}^{N_0+N-1} \left( \sum_{k \in \mathbb{Z}^{\star}}{C_\Gamma(\sigma_{n-kN-p})\ft{\mathscr{S}}(\sigma_{kN+p})}
\right)\\
& = \sum_{p=N_0}^{N_0+N-1}
{C_\Gamma(\sigma_{n-p})\ft{\mathscr{S}}(\sigma_{p})} + E_{\Gamma}(\sigma_{n}) \\
\end{split}
\label{eq:sampled_spectrum_irregular}
\end{equation}
where $E_{\Gamma}(\sigma_n)$ is the leakage in $\left[ N_0 \delta\sigma, \dots, (N_0 + N-1)\delta\sigma \right]$ from frequency samples outside this interval. Substituting $\mathscr{S}$ in equation~(\ref{eq:sampled_spectrum_irregular}) using equation~(\ref{eq:sampled_spectrum_regular}) and switching to the vectorial form yields:
\begin{equation}
\ft{{\cal S}}_{\Gamma} = N^{-1}\Delta x\,{\cal C}_{\Gamma}.\left[ \ft{{\cal S}}_{\mathrm{Reg}} -  {\cal E}_{\mathrm{Reg}} \right] + {\cal E}_\Gamma
\label{eq:vectorial_form}
\end{equation}
%





\noindent where ${\cal C}_\Gamma$ is the Toeplitz matrix:
\begin{equation}
{\cal C}_\Gamma=\left( \begin{array}{llll}
                       C_{\Gamma,{0}}  & C_{\Gamma,{-1}}    & \ldots & C_{\Gamma,{-N+1}} \\
                       C_{\Gamma,{1}}  & C_{\Gamma,{0}}     & \ldots & C_{\Gamma,{-N+2}} \\
                       \vdots & \vdots    & \ddots & \vdots     \\
                       C_{\Gamma,{N-1}} & C_{\Gamma,{N-2}} & \ldots & C_{\Gamma,{0}}
                     \end{array} \right) 
\end{equation}

\noindent whose coefficients are the values of the function $C_\Gamma$ sampled with steps of $\delta\sigma$:
\begin{equation}
C_{\Gamma,{k}}=C_\Gamma(\sigma_{k})={\Delta x}^{-1} \, \sum_{j=0}^{N-1}{e^{-2i\pi k\delta\sigma
x_{j}}}
\end{equation}
Hence the result and the proof of the lemma with $A=N^{-1}\Delta x\, {\cal C}_{\Gamma}$ and ${\cal R}={\cal E}_\Gamma - N^{-1}\Delta x\,{\cal C}_{\Gamma}.{\cal E}_{\mathrm{Reg}}$. \\ \\

This important lemma now leads to a method to retrieve the spectrum of the regularly sampled signal from the spectrum of the irregularly sampled signal. \\ \\
{\it{Theorem 1: inversion formula to recover the spectrum of the regularly sampled signal.}} The vector of the spectrum of the regularly sampled signal can be obtained from the spectrum of the irregular samples through the formula: $\ft{{\cal S}}_{\mathrm{Reg}}=N{\Delta x}^{-1}\,{\cal C}_{\Gamma}^{-1} \ft{{\cal S}}_{\Gamma}  + {\cal E}_{\mathrm{th}}$ with ${\cal E}_{\mathrm{th}} =  {\cal E}_{\mathrm{Reg}} - N{\Delta x}^{-1} \, {\cal C}_{\Gamma}^{-1}{\cal E}_\Gamma$. \\ \\
{\it{Proof}} \\
Given the result of the lemma, the proof consists in proving that the Toeplitz matrix is invertible. The matrix ${\cal C}_{\Gamma}$ can be factored into the product of a Vandermonde matrix ${\cal V}_{\Gamma}$ by its conjugate transpose:
\begin{equation}
{\cal C}_{\Gamma}={\Delta x}^{-1}\,{\cal V}_{\Gamma}.{\cal V}_{\Gamma}^{\star}
\end{equation}

\noindent with:
\begin{equation}
{\cal V}_{\Gamma}=\left( \begin{array}{llll} 1      & 1 & \ldots & 1 \\
e^{-2i\pi\delta\sigma x_{0}}  & e^{-2i\pi\delta\sigma x_{1}} & \ldots & e^{-2i\pi\delta\sigma
x_{N-1}} \\  
\vdots & \vdots    & \ddots & \vdots     \\ 
e^{-2i\pi(N-1)\delta\sigma x_{0}}& e^{-2i\pi (N-1)\delta\sigma x_{1}} & \ldots &
e^{-2i\pi (N-1)\delta\sigma x_{N-1}}
                     \end{array} \right) 
\end{equation}

\noindent whose determinant is equal to:
\begin{equation}
\mathrm{det}({\cal V}_{\Gamma}) = \prod_{N-1 \geqslant m>n\geqslant 0}{(e^{-2i\pi\delta\sigma x_{m}}-e^{-2i\pi\delta\sigma
x_{n}})}
\end{equation}

\noindent yielding for the determinant of ${\cal C}_{\Gamma}$:
\begin{equation}
\mathrm{det}({\cal C}_{\Gamma})=(\Delta x)^{-N}\prod_{N-1 \geqslant m>n\geqslant 0}
{|1-e^{2i\pi\delta\sigma(x_{m}-x_{n})}|^{2}}
\end{equation}
Since for all $N-1 \geqslant m>n\geqslant 0$, $x_{m}-x_{n} < \Delta x$ and $\delta\sigma(x_{m}-x_{n}) < 1$ by definition of $\Delta x$ and $\delta \sigma$, the determinant of ${\cal C}_{\Gamma}$ is non zero and equation~(\ref{eq:vectorial_form}) can be inverted yielding a method to recover the exact spectrum of the regularly sampled signal $S$ from the non-equispaced samples of $S_\Gamma$ through the relation:
\begin{equation}
\ft{{\cal S}}_{\mathrm{Reg}}=N{\Delta x}^{-1}\,{\cal C}_{\Gamma}^{-1} \ft{{\cal S}}_{\Gamma}  + {\cal E}_{\mathrm{th}}
\label{eq:general_solution}
\end{equation}
with ${\cal E}_{\mathrm{th}} =  {\cal E}_{\mathrm{Reg}} - N{\Delta x}^{-1} \, {\cal C}_{\Gamma}^{-1}{\cal E}_\Gamma$, thus proving theorem 1. \\ \\

It is to be noted that in the particular case of regular sampling ${\cal C}_{\mathrm{Reg}}=N{\Delta x}^{-1} I_N$ where $I_N$ is the identity matrix of size $N\times N$ and ${\cal E}_{\mathrm{th}} = 0$, and therefore $\ft{{\cal S}}_{\mathrm{Reg}}=\ft{{\cal S}}_{\Gamma}$.  
The possibility to recover the spectrum depends on the estimate of ${\cal E}_{\mathrm{th}}$. This is discussed in the next two sections in which the cases of periodic and non-periodic signals are addressed.
%
%
%
%
%
  \begin{figure}[t]
  \centering
\begin{tabular}{ll}
   \hspace{-0.6cm} \subfloat[Samples]{\label{fig:s_p_a}\includegraphics[height=\hauteur]{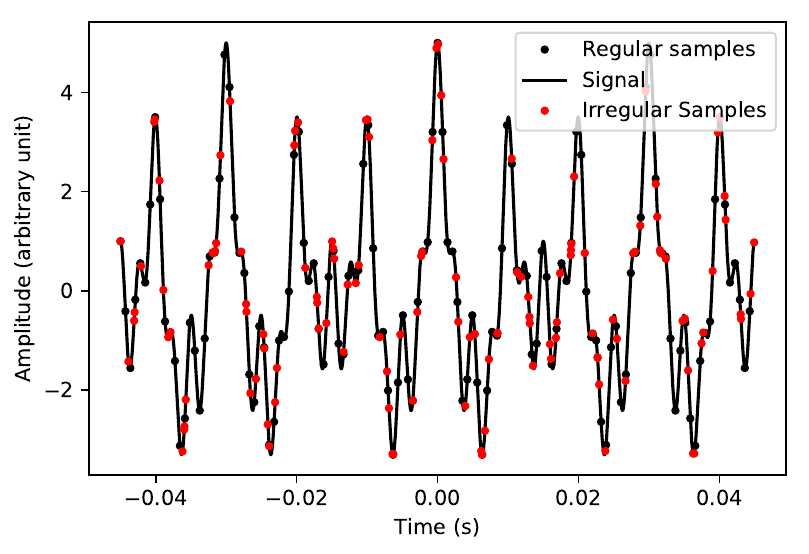}} & \hspace{-0.5cm}
        \subfloat[Spectra]{\label{fig:s_p_b}\includegraphics[height=\hauteur]{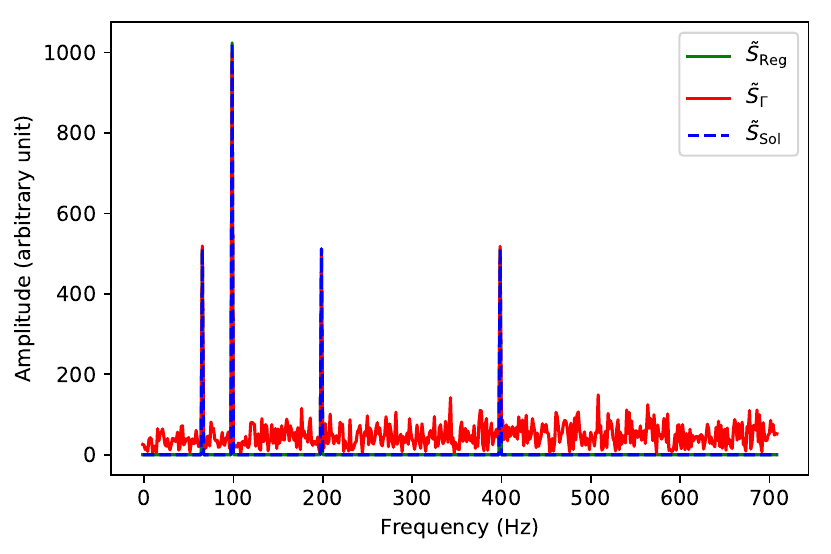}} \\    
      \hspace{-0.6cm} \subfloat[Sampling errors]{\label{fig:s_p_c}\includegraphics[height=\hauteur]{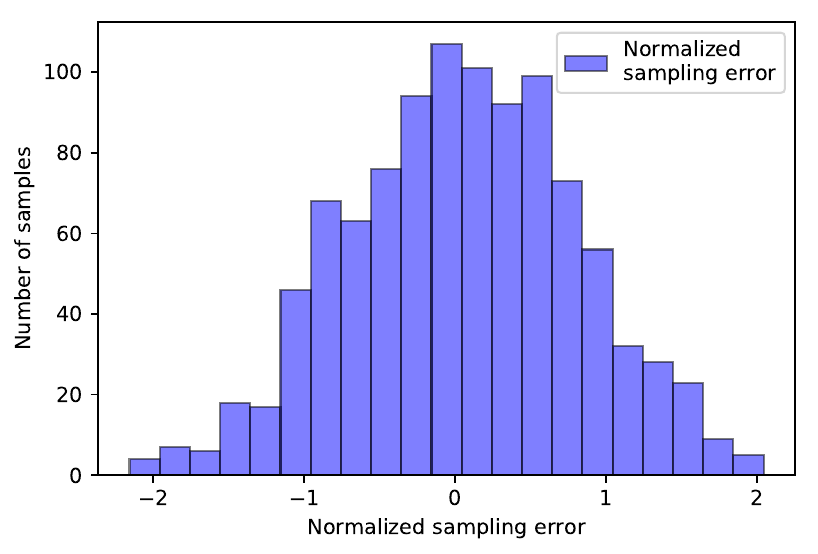}} & \hspace{-0.5cm}
      \subfloat[Log of spectra]{\label{fig:s_p_d}\includegraphics[height=\hauteur]{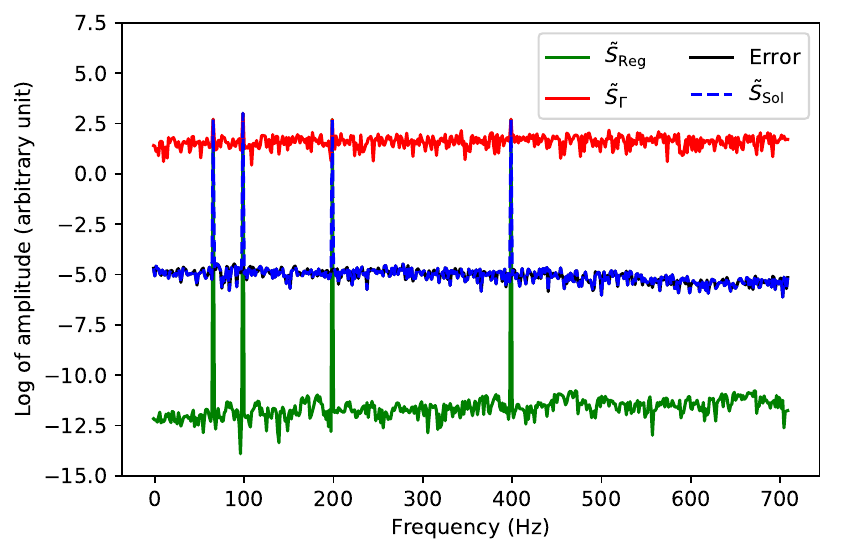}} \\
\end{tabular}
   \caption{Example of periodic band-limited signal: a mix of 4 frequencies (67, 100, 200 and 400 Hz) with respective amplitudes 1, 2, 1, 1. (a) Equispaced (black dots) and non-equispaced samples (red dots) plotted over theoretical signal; 1/8$^{\mathrm{th}}$ of the temporal window is displayed here. The amplitude of the sampling error in units of regular sampling interval is 2. (b) Amplitude of the spectrum $\ft{S}_\mathrm{Reg}$ of the signal with equispaced samples (green), amplitude of the direct DFT $\ft{S}_\Gamma$ of the non-equispaced samples (red) and amplitude of the spectrum $\ft{S}_\mathrm{Sol}$ reconstructed with the method of this paper (blue-dashed line). The condition number is $5.50\times10^{8}$. (c) Histogram of the normalized sampling errors. (d) Log of the spectrum. The black line is the amplitude of the difference $\ft{S}_{\mathrm{Sol}} - \ft{S}_{\mathrm{Reg}}$. }
         \label{Fig:samples_periodic}
   \end{figure}
\section{Computation of the spectrum for periodic band-limited signals}
\label{sec:spectrum_periodic}
In this section, $S$ is periodic and we assume that $\Delta x$ is a multiple of the period. \\ \\
%
{\it{Theorem 2: exact solution for the spectrum with non-equispaced data in the periodic band-limited case.}} In the particular case of a periodic band-limited signal sampled over a multiple of the period with the average sampling interval $\bar{\delta x} \leqslant \frac{1}{2f_\mathrm{max}}$, the spectrum of the regularly sampled signal is exactely deduced from the spectrum of the irregularly sampled signal through the inverse of a Toeplitz matrix:
\begin{equation}
\ft{{\cal S}}_{\mathrm{Sol}}=\ft{{\cal S}}_{\mathrm{Reg}}=N{\Delta x}^{-1}\,{\cal C}_{\Gamma}^{-1} \ft{{\cal S}}_{\Gamma}
\label{eq:exact_solution_periodic}
\end{equation}
\label{th:spectrum_periodic}
%
{\it{Proof}} \\
Given the result of theorem 1, proving theorem 2 is equivalent to showing that ${\cal E}_{\mathrm{th}}=0$ in the particular case of a periodic signal sampled over a multiple of the period. In this case, the samples of the discrete Fourier transform of $\mathscr{S}$ are proportional to the samples of the discrete Fourier transform of $S$ as there is no frequency leakage~\cite{Brigham1988}. $\mathscr{S}$ is therefore also band limited and its spectrum is equal to 0 outside the interval $[-f_\mathrm{max},+f_\mathrm{max}]$, implying that $\ft{\mathscr{S}}(\sigma_{k})$ is zero for $k$ outside the interval $[N_0, N_0+N-1]$. 
As a consequence, the sum in equation~(\ref{eq:sampled_spectrum}) has a limited number of terms and writes:
\begin{equation}
\ft{S}_{\Gamma}(\sigma_{n})=\sum_{k=N_0}^{N_0+N-1}{C_\Gamma(\sigma_{n-k})\ft{\mathscr{S}}(\sigma_{k})}
\label{eq:sampled_irregular_spectrum_periodic}
\end{equation}
with $E_{\Gamma}(\sigma_{n}) = 0$ for all $n$ between $N_0$ and $N_0+N-1$. 
Likewise, in the regular sampling case and for $n$ in the same interval, equation~(\ref{eq:sampled_spectrum_regular}) becomes:
\begin{equation}
\ft{S}_{\mathrm{Reg}}(\sigma_{n}) = N{\Delta x}^{-1}{\ft{\mathscr{S}}(\sigma_{n})}
\label{eq:sampled_spectrum_regular_periodic}
\end{equation}
with $E_{\mathrm{Reg}}(\sigma_{n}) = 0$. For all $n$ in the interval above we therefore have:
\begin{equation}
\ft{S}_{\Gamma}(\sigma_{n})=N^{-1} \Delta x \sum_{k=N_0}^{N_0+N-1}{C_\Gamma(\sigma_{n-k})\ft{S}_{\mathrm{Reg}}(\sigma_{k})}
\end{equation}
and the column vector $\cal{E}_\mathrm{th}$ is equal to 0, hence the result of theorem 2. \\ \\
%
%
%
%
{\it{Remark 3: }} similar or identical results were found for the reconstruction of periodic functions from non-uniform samples \cite{Cheng2020,Feichtinger1995}. \\ \\

Fig.~\ref{Fig:samples_periodic} shows an example of a spectrum of a periodic signal $\ft{S}_{\mathrm{Sol}}$ reconstructed with the method presented in this paper. The theoretical function is the sum of four cosine functions whose frequencies are harmonics of $\frac{100}{3}$\,Hz:
\begin{equation}
\label{eq:periodic}
S(t)= \cos\left(\frac{2\pi t}{0.0025\,\mathrm{s}}\right)+\cos\left(\frac{2\pi t}{0.005\,\mathrm{s}}\right)+2\cos\left(\frac{2\pi t}{0.01\,\mathrm{s}}\right)+\cos\left(\frac{2\pi t}{0.015\,\mathrm{s}}\right)
\end{equation}
The comparison with the spectrum of the regular grid samples  $\ft{S}_\mathrm{Reg}$ and with the result of the direct Discrete Fourier Transform $\ft{S}_\Gamma$ are also plotted. 1024 samples have been used between -0.72 and +0.72\,s (only the central 0.1\,s is plotted). The irregular samples have been derived from the set of regular samples by adding randomly distributed sampling errors between -2$\delta x$ and +2$\delta x$, with $\delta x$ the sampling interval before the errors were added. The histogram of the normalized sampling errors is shown with sampling errors in units of $\delta x$.
Other than the method, one finding of this paper is that the condition number of the Toeplitz matrix plays a critical role in the accuracy of the result (see \ref{sec:accuracy}). For this particular example the condition number is 5.50$\times 10^{8}$. The condition number is dependent on the realization of the sampling errors and is not constant for a given statistics of the errors. The condition number therefore varies for a fixed distribution law. In the example presented here, the dynamic range of the reconstructed spectrum is a few $10^7$ compared to a few 10 with the direct Discrete Fourier Transform, hence an improvement of the order of $10^6$ in dynamic range with the method proposed in this paper.
   
 \begin{figure}[t]
 \centering
\begin{tabular}{ll}
     \hspace{-0.6cm} \subfloat[Samples]{\label{fig:s_i_a}\includegraphics[height=\hauteur]{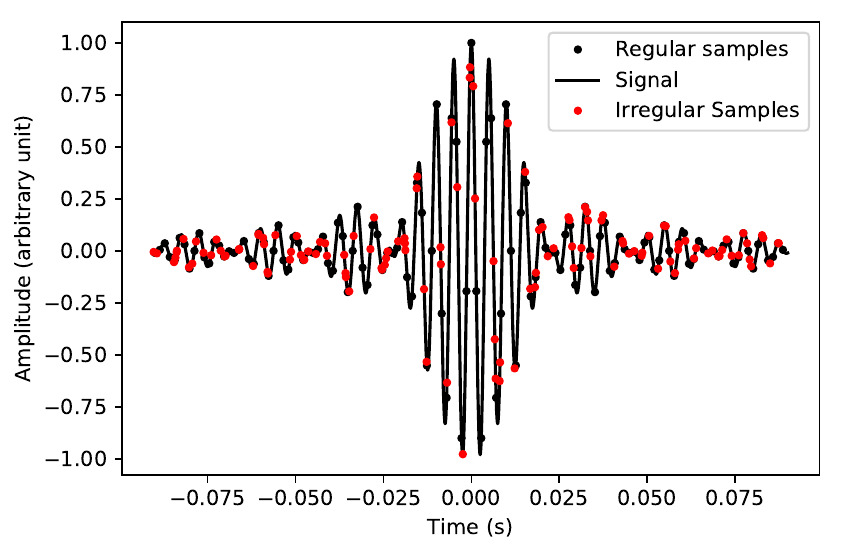}} & \hspace{-0.5cm}
        \subfloat[Spectra]{\label{fig:s_i_b}\includegraphics[height=\hauteur]{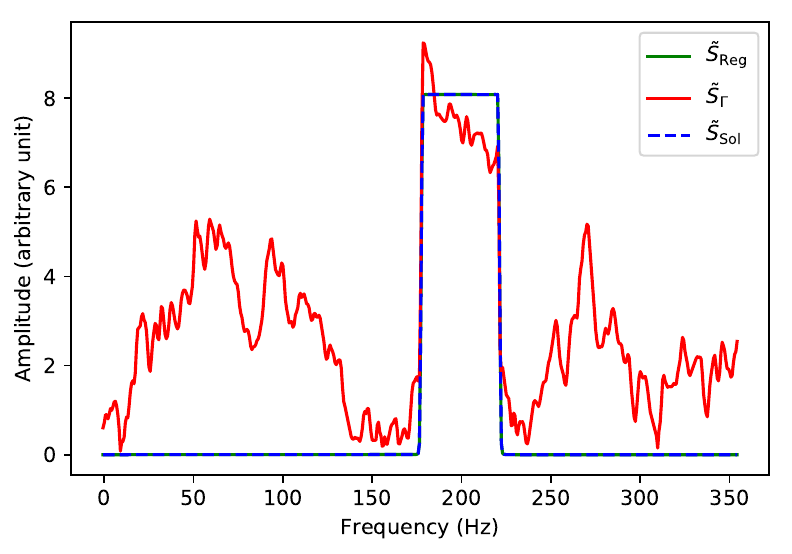}} \\    
       \hspace{-0.6cm} \subfloat[Sampling errors]{\label{fig:s_i_c}\includegraphics[height=\hauteur]{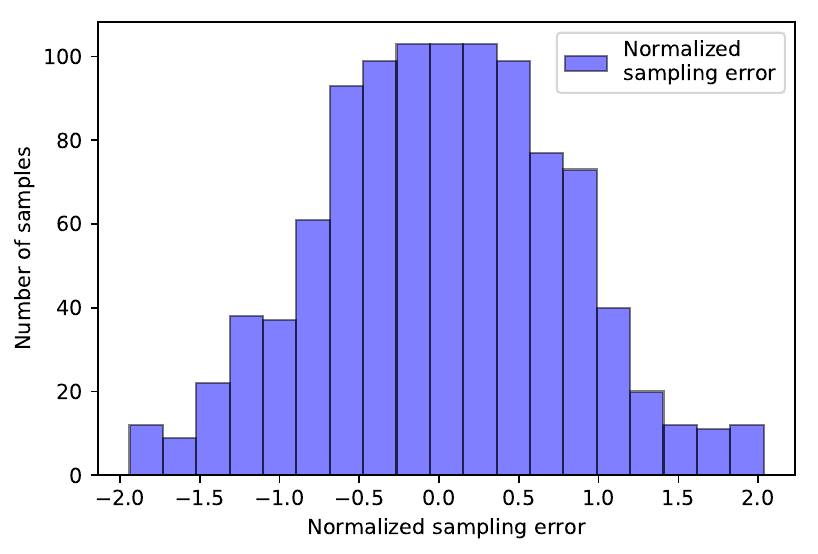}} & \hspace{-0.5cm}
      \subfloat[Log of spectra]{\label{fig:s_i_d}\includegraphics[height=\hauteur]{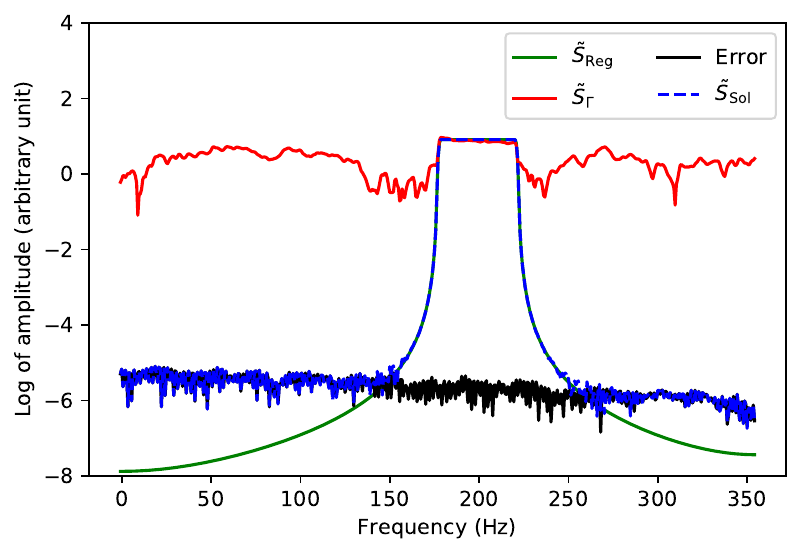}} \\
\end{tabular}
   \caption{Example of non-periodic band-limited function: a fringe pattern with 200\,Hz temporal frequency and width 44\,Hz. (a) Equispaced (black dots) and non-equispaced samples (red dots) plotted over theoretical signal. The amplitude of the sampling error in units of regular sampling interval is 2. (b) Amplitude of the spectrum $\ft{S}_\mathrm{Reg}$ of the signal with equispaced samples (green), amplitude of the direct DFT $\ft{S}_\Gamma$ of the non-equispaced samples (red) and amplitude of the spectrum $\ft{S}_\mathrm{Sol}$ reconstructed with the method of this paper (blue-dashed line). The condition number is $4.43\times10^{10}$. (c) Histogram of the normalized sampling errors. (d) Log of the spectrum. The black line is the amplitude of the difference $\ft{S}_{\mathrm{Sol}} - \ft{S}_{\mathrm{Reg}}$.}
         \label{Fig:samples_interferogram}
   \end{figure}
%
%
%
%
\section{Computation of the spectrum for non-periodic band-limited signals}
\label{sec:spectrum_non-periodic}
If $S$ is band-limited but not periodic then $\mathscr{S}$ is not band-limited, whatever the length of the sampling range. The exact spectrum cannot be deduced from the spectrum of the non-equispaced samples as shown in equation~(\ref{eq:general_solution}). An approximated solution $\ft{S}_{\mathrm{Sol}}$ can be obtained  ignoring the error term and directly using the result for the periodic case of theorem~2. Doing so, a theoretical error of $E_{\mathrm{th}}$ is made on the spectral samples.
The magnitude of the error $E_{\mathrm{th}}$ depends on the amount of frequency leakage which depends on the properties of the truncated function. Leakage can be reduced by apodizing the raw signal \cite{Brigham1988}. Some functions are naturally apodized and have very low leakage. 
An example is shown in Fig.~\ref{Fig:samples_interferogram} which is equivalent to Fig.~\ref{Fig:samples_periodic} but with the non-periodic function: 
\begin{equation}
\label{eq:interf}
S(t) = \cos\left(\frac{2\pi t}{0,005\,\mathrm{s}} \right)\times\sinc \left(\frac{t}{0,023\,\mathrm{s}}\right)
\end{equation}
where $\sinc{(x)}=\sin(\pi x)/(\pi x)$. The spectrum of the untruncated signal $S$ is a pure door function of width 44\,Hz centered on 200\,Hz. It is the intensity of an interferogram recorded with a Michelson interferometer with a scanning mirror on a source with a finite and flat spectral bandwidth. The effect of leakage is visible in the bottom right panel of Fig.~\ref{Fig:samples_interferogram} with wings showing up in the logarithmic scale. Leakage has been reduced in this example by multiplying the sequence by a Hann window~\cite{Harris1978}. The use of the Hann window also reduces the effect of ringing due to the sharp edges in the spectral domain in this particular case. The difference between the spectrum of the truncated regularly sampled signal $\ft{S}_\mathrm{Reg}$ and the spectrum reconstructed with the method of this paper $\ft{S}_{\mathrm{Sol}}$ is shown in black. Contrary to the example of the periodic function, the direct Discrete Fourier Transform of the irregular samples $\ft{S}_\Gamma$ gives a very poor result with a dynamic range of order 2 while the dynamic range in $\ft{S}_{\mathrm{Sol}}$ is of the order of a few $10^6$. The same statistics of the sampling error has been used as for the periodic function in Section~\ref{sec:spectrum_periodic} with a condition number for the particular realization presented here of $4.43\times10^{10}$. The example shows how accurate the method can be, even for non-periodic band-limited functions. We have tried another non-periodic function and applied the method to the Ricker wavelet as in \cite{Gulati2010}:
\begin{equation}
\label{ricker}
S(t) = \exp \left(- \left(\frac{\pi t}{1\,\mathrm{s}} \right)^2 \right) * \left( 1-2 \left(  \frac{\pi t}{1\,\mathrm{s}} \right)^2   \right)
\end{equation}
The Ricker function is naturally apodized because of the Gaussian envelope and no filtering with the Hann window or an equivalent was necessary. The results obtained with this particular function are discussed in Section~\ref{sec:accuracy}.
%
%
%
%
\section{Complexity and accuracy of the algorithm}
\label{sec:accuracy}
\begin{figure}[h]
\centering
\includegraphics[width=10cm]{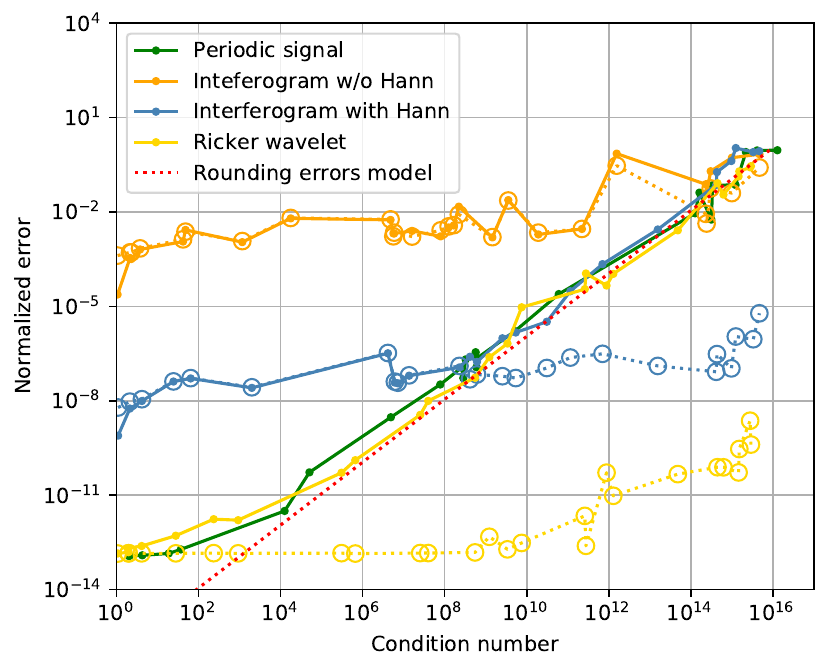}
   \caption{Plot of the normalized errors on the reconstructed spectra vs. the condition number of the Toeplitz matrix for various realizations of sampling grids for the periodic signal of Section~\ref{sec:spectrum_periodic} (green) and the non-periodic signals of Section~\ref{sec:spectrum_non-periodic} (interferogram in orange, interferogram with Hann windowing in blue and Ricker wavelet in yellow). The errors have been normalized by the maximum of the spectrum amplitude to compare the impact of relative machine rounding errors on the reconstructed spectra. The continuous lines are the errors $\|\ft{{\cal S}}_{\mathrm{Sol}} - \ft{{\cal S}}_{\mathrm{Reg}}\|$ normalized by the maximum of  $|\ft{{S}}_{\mathrm{Reg}}|$. The red-dotted line is a model of the trend of the error due to machine precision as explained in Section~\ref{sec:accuracy}. The dotted lines with open circles are the theoretical errors ${\cal E}_{\mathrm{th}}$ of equation~(\ref{eq:general_solution}). Note that the computation of the theoretical errors is also subject to machine precision beyond condition number $10^{15}$.}
         \label{Fig:errors_condition}
   \end{figure}
The computation of the spectra has been coded in \texttt{python}. The \texttt{solve\_toeplitz} routine from the linear algebra \texttt{scipy.linalg} has been used to solve for $\ft{{\cal S}}_{\mathrm{Sol}}$ in equation~(\ref{eq:exact_solution_periodic}). It is based on the Levinson-Durbin recursion and runs in $\mathcal{O}(N^2)$~\cite{Levinson1946}. In the case of positive-definite Toeplitz matrixes as here, the algorithm can be optimized to run in $\mathcal{O}(N\log^{2}N)$ \cite{Ammar1996} \msgguy{and even in $\mathcal{O}(N\log N)$ \cite{Chan1989} with the the preconditioned conjugate gradient method.}\\

\noindent The errors on the reconstructed spectra $\|\ft{{\cal S}}_{\mathrm{Sol}} - \ft{{\cal S}}_{\mathrm{Reg}}\|$ normalized by the maximum of  $|\ft{{\cal S}}_{\mathrm{Reg}}|$, where $\|x\|$ is the 2-norm or Euclidian norm of $x$, are plotted in Fig.~\ref{Fig:errors_condition}. It is remarkable that the various signals lead to similar behaviors. First of all, above a certain threshold, the errors scale with the condition number. Below this threshold, the errors tend to stick to a minimum floor. This can be understood from equation~(\ref{eq:general_solution}) where the first term is the exact solution for periodic signals and the approximate solution for non-periodic signals while the second term ${\cal E}_{\mathrm{th}}$ is the theoretical error for non-periodic functions, equal to zero in the periodic case. At low condition number, the second term dominates while at high condition number, the error on the first term  due to machine rounding errors prevails.

\noindent For the periodic function, and in absence of rounding errors, the solution is exact as  ${\cal E}_{\mathrm{th}}=0$ and  $\ft{{\cal S}}_{\mathrm{Sol}}=\ft{{\cal S}}_{\mathrm{Reg}}=N{\Delta x}^{-1}\,{\cal C}_{\Gamma}^{-1} \ft{{\cal S}}_{\Gamma}$. Limited machine precision however introduces rounding errors that propagate in the computation of the solution even at low condition number with errors in computing ${\cal C}_{\Gamma}^{-1}$ and $\ft{{\cal S}}_{\Gamma}$. This leads to the minimum relative error of order 1e-13 in Fig.~\ref{Fig:errors_condition} at very low condition number. The theoretical error on the Ricker wavelet reconstructed spectrum (yellow dashed line and open yellow circles in  Fig.~\ref{Fig:errors_condition}) being very close to this minimum error, both Ricker and periodic function error curves have the same behavior at low condition number and the effect of the theoretical error is not visible on the graph in this particular case.
%
%
\noindent At high condition number, the error is dominated by the ill-conditioned nature of the inversion of ${\cal C}_{\Gamma}$. An upper bound is obtained by applying the definition of the condition number:
\begin{equation}
\frac{\|\Delta\ft{{\cal S}}_{\mathrm{Sol}}\|}{\|\ft{{\cal S}}_{\mathrm{Sol}}\|} \leqslant \kappa({\cal C}_{\Gamma}) \frac{\|\Delta\ft{{\cal S}}_{\Gamma}\|}{\|\ft{{\cal S}}_{\Gamma}\|}
\end{equation}
From computation, we have numerically checked that $\|\ft{{\cal S}}_{\mathrm{Sol}}\| \approx {\|\ft{{\cal S}}_{\Gamma}\|}$ which implies that:
\begin{equation}
\|\Delta\ft{{\cal S}}_{\mathrm{Sol}}\| \lessapprox \kappa({\cal C}_{\Gamma}) \|\Delta\ft{{\cal S}}_{\Gamma}\|
\end{equation}
An accurate assessment of the propagation of the rounding errors goes beyond the scope of this paper. The error analysis is limited here to get a trend. The machine error is assigned to macroscopic ensembles of operations, not reaching the level of each individual operations of the computation of the spectra. The numerical error on $\ft{{\cal S}}_{\Gamma}$ is equal to $\|\Delta\ft{{\cal S}}_{\Gamma}\|= \sqrt{\sum_{i=0}^{N-1}{|\Delta\ft{{\cal S}}_{\Gamma,i}} |^2}$. Each $\ft{{\cal S}}_{\Gamma,i}$ is the sum of $N$ terms leading to a real and an imaginary part with associated errors $\Delta\mathrm{Re}(\ft{{\cal S}}_{\Gamma,i})$ and $\Delta\mathrm{Im}(\ft{{\cal S}}_{\Gamma,i})$, respectively. Each sample $S(x_k)$ contributes a rounding error proportional to the relative error $\epsilon_\mathrm{m}$ in the computation of the Discrete Fourier Transform $\ft{{\cal S}}_{\Gamma,i}$. In the case of the periodic signal, the samples $S(x_k)$ can be considered of equal strength of order 1 (it is not the case for the interferogram and the Ricker wavelet which are localized close to $t=0$\,s with respectively $\simeq 128$ and 40 samples effectively building the error). In the case of the periodic function the variance of rounding errors is therefore:
\begin{equation}
\mathrm{var}(\Delta\mathrm{Re}(\ft{{\cal S}}_{\Gamma,i}))=\mathrm{var}(\Delta\mathrm{Im}(\ft{{\cal S}}_{\Gamma,i}))=\mathcal{O}(N)\epsilon^2_\mathrm{m}
\end{equation}
 Assuming that the statistics of the rounding errors are Gaussian, then $\|\Delta\ft{{\cal S}}_{\mathrm{Sol}}\|^2$ is a $\chi^2$ with $2N$ degrees of freedom with average value $2N\times \mathcal{O}(N)\epsilon^2_\mathrm{m}$. Since $\langle \|\Delta\ft{{\cal S}}_\Gamma\| \rangle \leq \sqrt{\langle  \|\Delta\ft{{\cal S}}_\Gamma\|^2 \rangle}$ we therefore get:
 %
%
\begin{equation}
\label{eq:error_high}
\|\Delta\ft{{\cal S}}_{\mathrm{Sol}}\|  \lessapprox \kappa({\cal C}_{\Gamma}) \mathcal{O}(N) \times \epsilon_\mathrm{m}
\end{equation}
In the case of the periodic function, the normalization of the spectrum leads to dividing by the number of samples $N$. Hence the approximation of the upper bound for the relative error which is proportional to $ 2^{-53}\times\kappa({\cal C}_{\Gamma})$ for double-precision floats. This model is plotted on \msgguy{Fig.}~\ref{Fig:errors_condition} as a red-dotted line and nicely reproduces the trend  of the relative error due to machine precision at high condition number. \\ 
For non-periodic functions, the minimum error is set by the theoretical error ${\cal E}_{\mathrm{th}}$ of equation~(\ref{eq:general_solution}). The theoretical error is plotted with open circles and dotted lines on Fig.~\ref{Fig:errors_condition}. It matches quite well the difference between the reconstructed spectrum and the spectrum directly obtained with regular samples for the lowest range of condition numbers. The transition to the error set by the machine precision clearly depends on the amount of apodization in the function. The error for the interferogram without Hann filtering gets limited by machine precision only beyond a condition number of $10^{14}$ while it is as early as $10^{9}$ with a Hann window, showing the interest of this filtering to increase the accuracy of the method. As seen above, the Ricker wavelet is naturally more apodized due to the Gaussian component with a much lower theoretical error and a transition at $10^{3}$ if rounding errors would not dominate with double-precision floats at low condition number. The computation of the theoretical error ${\cal E}_{\mathrm{th}}$  is also prone to rounding errors which explains why the trend beyond $10^{15}$  gets parallel to the red dotted line in the Log-Log graph. All cut-off condition numbers would be shifted to higher values with increased machine precision (by $\times 2048$ for extended precision and by $\times 2^{60}$ for quadruple precision). \msgguy{In} the latter case, none of the examples shown here would be affected by machine precision with condition numbers limited to $10^{16}$ and Hann filtering would allow to reach dynamic ranges as high as $10^7$ with the method. \msgguy{An alternative method to the direct inversion of the Toeplitz system is worth investigating though to improve the accuracy of the method but goes beyond the scope of this paper.}
\section{Comparison with other methods}
\label{section:comparison}
\noindent\msgguy{The performance of the method presented in this paper can be compared to the other methods cited in the introduction in Section~\ref{section:introduction}. The iterative methods allow to reconstruct the samples/spectra with accuracies comparable to what is achieved in this paper. The multichannel method requires additional data like samples of the derivative of the signal to achieve accurate reconstructions. Some iterative methods are designed for some specific cases of non-regular sampling like in \cite{Strohmer2006} for periodic nonuniform sampling. In practice, only the iterative methods based or derived from the iterative methods by means of the frame operators of \cite{Duffin1952} that reconstruct the signal from the samples of the signal only and are applied to the general case of nonuniform sampling can be compared to the results of this paper. Comparisons of iteratives methods are discussed for example in \cite{Feichtinger1995, feichtinger2021}. All methods require the maximum interval between two samples to be less than the Nyquist interval to get an estimate of the rate of convergence and therefore to ensure convergence. It is not the case of the method of this paper for which the consequence of gaps larger than the Nyquist interval is to increase the condition number of the Toeplitz matrix and reduce precision but not necessarily to prevent reconstruction. The adaptive weights method of \cite{feichtinger2021} and the improved versions with acceleration by the conjugated gradients method allow to reach relative accuracies of $10^{-13}$ in ~20-40 iterations depending on the method~\cite{Feichtinger1995}. It would be necessary to compute the condition number of the Toeplitz matrix for the particular example of the sampled signal in this paper to make a direct comparison. However, one may take the two cases of Fig.~\ref{Fig:samples_periodic} and Fig.~\ref{Fig:samples_interferogram} of the present paper for which accuracies of order $10^{-7}$ and $10^{-5}$ are respectively achieved with condition numbers of $5.5\times 10^8$ and $4.4 \times 10^{10}$. Such accuracies would require of the order of 12 and 8 iterations with the best algorithm (PACT, see below) according to Fig. 3 of  \cite{Feichtinger1995}. Although the convergence rate of the iterative algorithms cannot be predicted when the maximum gap exceeds the Nyquist interval, \cite{Feichtinger1995} show that the adaptive weight method with congugated gradients (ACT) and the version upgraded with the use of a preconditioner (PACT) can be successful to reconstruct signals in this case though, at least in the instance presented in Fig.~5 of their paper. A similar case could be approached with the statistics of the samples for the interferogram of Section~\ref{sec:spectrum_non-periodic} with about 10\% of the gaps between samples being larger than the Nyquist interval. For this particular realization of the samples, this yielded a condition number of $7.3\times 10^{12}$ and a relative error for the reconstructed spectrum of the order of $10^{-3}$. According to Fig.~5 of their paper, such precision requires more than 230 iterations with the ACT method and is reached in 130 iterations with the PACT method. \\ 
As far as complexity is concerned, as the iterative algorithms also require to solve a Toeplitz system, their complexity is therefore the same as in this paper except that they need several iterations to converge while a single step is required here. But their accuracy is not limited by the theoretical error ${\cal E}_{\mathrm{th}}$ of equation~(\ref{eq:general_solution}) and signals can be {\it a priori} reconstructed with an accuracy limited by machine precision only. One may therefore contemplate to combine the two methods by performing the first iteration with the method presented in this paper and iterating from there on with the iterative methods. Depending on the characteristics of the signal and on the sampling sequence, it may reduce the number of steps by a significant factor. A more comprehensive comparison of the methods is necessary at this stage to reach further conclusions.}
%
%
%
%
%
\section{Generalization to higher dimensions}
\label{sec:dimension} The method can be easily generalized to any dimension $d \geqslant 1$. Sampling points are now vectors $\vx=[x_1, \cdots, x_d]$ of ${\Bbb{R}}^{d}$ (vectors are in bold face).
%
%
The grid of samples writes $\Gamma=\{ {\vx}_{0},\ldots,{\vx}_{N^{d}-1} \}$ where we choose to sort the points in lexicographic order. For the sake of simplicity, we have assumed that the number of points is $N^{d}$ with $N$ points in average per dimension. The sampled function writes:
\begin{equation}
S_\Gamma({\vx})=\sum_{i=0}^{N^{d}-1}{S(\vx_{i})\delta(\vx-\vx_{i})}
\end{equation} 
and the spectrum:
\begin{equation}
\ft{S}_\Gamma(\vsig)=\sum_{i=0}^{N^{d}-1}{S(\vx_{i})e^{-2i\pi\vsig.\vx_{i}}}
\end{equation}
where $\vsig=[\sigma_1,\cdots,\sigma_d]$ is the $d$-frequency conjugated to $\vx$. The average sampling interval is also a $d$-vector of components:
\begin{equation}
\overline{{\delta x}_{i}}=\frac{ \max_{m,n}|{ x_{i,m}-x_{i,n} }| }{ N-1 }
\end{equation}
with the sampling width $\Delta \vx = [N \overline{{\delta x}_{1}},\cdots,N \overline{{\delta x}_{d}}]$. As in dimension 1, the spectral sampling interval is:  $\vdel\vsig=[({\Delta x}_{1})^{-1},\cdots,({\Delta x}_{d})^{-1}]$. With these notations, all relations derived in dimension 1 can be generalized to any dimension $d$ by replacing the product of $x$ by the associated frequency $\sigma$ by a scalar product. The grid function becomes:
\begin{equation}
G_\Gamma(\vx)=\sum_{i=0}^{N^{d}-1}{\left[ \sum_{{\bf n} \in \mathbb{Z}^d}{\delta (\vx-\vx_{i}-{\bf n}\had\Delta \vx)} \right]}
\end{equation}
where ${\bf n}$ is a $d$-index of components $[n_{1},n_{2},...,n_{d}]$ and $\bf{A}\had\bf{B}$ is the Hadamard product of vectors $\bf{A}$ and $\bf{B}$. The Fourier transform of the grid distribution is equal to:
\begin{equation}
{\ft G}_\Gamma(\vsig)=({\Delta x}_{1}\dots{\Delta x}_{d})^{-1} \left[ \sum_{i=0}^{N^{d}-1}{e^{-2i\pi \vsig
.\vx_{i}}}\right] 
                      \diracComb_{\delta \sigma_{1} \dots \delta \sigma_{d}} (\sigma_{1},\dots,\sigma_{d})
\end{equation}
with $\diracComb_{\delta \sigma_{1} \dots \delta \sigma_{d}}$ the $d$-dimension Dirac Comb of steps $\delta \sigma_{1} \dots \delta \sigma_{d}$. Hence the $C_\Gamma$ function:
\begin{equation}
C_\Gamma(\vsig)=({\Delta x}_{1}\dots{\Delta x}_{d})^{-1} \left[ \sum_{i=0}^{N^{d}-1}{e^{-2i\pi \vsig .
\vx_{i}}}\right] 
\end{equation}
and the spectrum of the grid distribution:
\begin{equation}
\ft{G}_\Gamma(\vsig)=
\sum_{n_{1},n_{2},...,n_{d}=-\infty}^{+\infty}{C_\Gamma({\bf n}\had{\mbox{\boldmath
$\delta\sigma$}})\delta(\vsig-{\bf n}\had{\mbox{\boldmath $\delta\sigma$}})} 
\end{equation}
The sampled signal writes:
\begin{equation}
S_\Gamma(\vx)=S_\Gamma(\vx).\Pi(\vx)G_\Gamma(\vx)={\mathscr{S}}(\vx)G_\Gamma(\vx)
\end{equation}
The sampled spectrum writes:
\begin{equation}
\ft{S}_\Gamma({\bf n}\had{\mbox{\boldmath $\delta\sigma$}})=\ft{S}_\Gamma({\mbox{\boldmath $\sigma_{\bf n}$}})=\sum_{{\bf k} \in \mathbb{Z}^{d}}{C_\Gamma{\mbox{\boldmath
$(\sigma_{\bf n-k}$}})\ft{\mathscr{S}}({\mbox{\boldmath $\sigma_{\bf k}$}})} 
\end{equation}
where ${\bf k}$ is a $d$-index. It is assumed that the signal is at least Nyquist sampled, ie $\overline{{\delta x}_{i}} \leqslant \frac{1}{2 f_{\mathrm{max},i}}$ for all $i=1,\dots,d$. The components of the sampled frequencies are in the range $\left[ N_0 \delta\sigma, \dots, (N_0 + N-1)\delta\sigma \right]$.
%
%
%
 Like in equations~(\ref{eq:sampled_spectrum_regular}) and (\ref{eq:sampled_spectrum_irregular}), the above expression can be split in two parts, yielding for the spectrum on a regular grid:
\begin{equation}
\begin{split}
\ft{S}_{\mathrm{Reg}}(\sigma_{\bf n}) & = N^d (\Delta x_1 \cdots \Delta x_d)^{-1}{\ft{\mathscr{S}}(\vsig_{\bf n})} + N^d (\Delta x_1 \cdots \Delta x_d)^{-1}\sum_{{\bf k} \in \mathbb{Z}^{\star d}}{\ft{\mathscr{S}}(\vsig_{N{\bf k}+{\bf n}})}\\
& = N^d (\Delta x_1 \cdots \Delta x_d)^{-1}{\ft{\mathscr{S}}(\vsig_{\bf n})} + E_{\mathrm{Reg}}(\vsig_{\bf n})
\end{split}
\label{eq:sampled_spectrum_regular_d}
\end{equation}
and for an irregular grid:
\begin{equation}
\begin{split}
\ft{S}_{\Gamma}(\sigma_{\bf n}) & = \sum_{p_1 \cdots p_d=N_0}^{N_0+N-1}
{C_\Gamma(\vsig_{{\bf n}-{\bf p}})\ft{\mathscr{S}}(\vsig_{\bf p})} + \sum_{p_1 \cdots p_d=N_0}^{N_0+N-1} \;\; \sum_{{\bf k} \in \mathbb{Z}^{\star d}}{C_\Gamma(\vsig_{{\bf n}-N{\bf k} -{\bf p}})\ft{\mathscr{S}}(\vsig_{N{\bf k}+{\bf p}})}\\
& = \sum_{p_1 \cdots p_d=N_0}^{N_0+N-1}
{C_\Gamma(\vsig_{{\bf n}-{\bf p}})\ft{\mathscr{S}}(\vsig_{\bf p})} + E_{\Gamma}(\vsig_{\bf n}) \\
\end{split}
\label{eq:sampled_spectrum_irregular_d}
\end{equation}
One recognizes in the above formula a linear relation between the regularly sampled spectrum, the spectrum of the samples on the grid $\Gamma$ and an error term. One can therefore write the equivalent of equation~(\ref{eq:general_solution}) and define the solution spectrum as:
\begin{equation}
\label{eq:sol_dim_d}
\ft{{\cal S}}_{\mathrm{Sol}}=N^d (\Delta x_1 \cdots \Delta x_d)^{-1}\,{\cal C}_{\Gamma}^{-1} \ft{{\cal S}}_{\Gamma}  \end{equation}
the theoretical error being ${\cal E}_{\mathrm{th}} =  {\cal E}_{\mathrm{Reg}} - N^d (\Delta x_1 \cdots \Delta x_d)^{-1} \, {\cal C}_{\Gamma}^{-1}{\cal E}_\Gamma$.  Sorting the $d$-indices in lexicographic order, the Toeplitz matrix ${\cal C}_\Gamma$ writes:
\begin{equation}
{\cal C}_\Gamma=\left( \begin{array}{llll}
                       C_{\Gamma,{[0,\cdots,0]}}  & C_{\Gamma,{[-1,0,\cdots,0]}}    & \ldots & C_{\Gamma,{[-N+1,\cdots,-N+1]}} \\
                       C_{\Gamma,{1,[0,\cdots,0]}}  & C_{\Gamma,{[0,\cdots,0]}}     & \ldots & C_{\Gamma,{[-N+2, -N+1,\cdots,-N+1]}} \\
                       \vdots & \vdots    & \ddots & \vdots     \\
                       C_{\Gamma,{[N-1,\cdots,N-1]}} & C_{\Gamma,{[N-2, N-1,\cdots,N-1]}} & \ldots & C_{\Gamma,{[0,\cdots,0]}}
                     \end{array} \right) 
\end{equation}
In the particular case of regular sampling ${\cal C}_{\mathrm{Reg}}=N^d(\Delta x_1\cdots \Delta x_d)^{-1}I_N^d$ where $I_N^d$ is the identity matrix of size $N^d\times N^d$ and ${\cal E}_{\mathrm{th}} = 0$. Similarly to the 1D case, the solution is exact for periodic functions sampled over an integer number of periods and the theoretical error ${\cal E}_{\mathrm{th}} = 0$. All conclusions reached in the 1D case in the previous sections are  applicable to any dimension $d$. 
%
%
%
%
\section{Conclusion}
\label{sec:conclusion}
A method has been presented in this paper to compute the discrete spectrum of a signal from $N$ non-equispaced samples. In the case of a periodic function sampled over a multiple of the period and respecting the Nyquist criterion, the discrete spectrum is linearly linked to the direct DFT of the samples by an $N \times N$ Toeplitz matrix thereby providing an exact solution to the problem of computing the spectrum. The Toeplitz matrix is positive-definite and can be inverted in \msgguy{$\mathcal{O}(N\log N)$} operations. This same method can be applied to non-periodic signals with an excellent approximation whose theoretical error is derived. The accuracy can be improved by a factor $10^3$ by filtering the signal with a Hann window in the example chosen in the paper. The inversion of the Toeplitz matrix is ill-conditioned leading to increasing errors  with increasing irregularity of the samples. It has been shown that the error on the reconstructed spectrum is limited by machine precision beyond a certain threshold of condition number which depends on the nature of the signal. Spectra can be reconstructed with dynamic ranges as high as $10^{13}$ with double-float precision. \msgguy{The method has been compared to iterative methods. The complexity is similar albeit with a single iteration.} \corguy{High up to very high relative accuracies for the reconstructed spectrum can be reached in a single step,} \msgguy{but iterative methods may reach higher accuracies as they are not limited by the theoretical error derived in this paper for non-periodic signals. If higher precisions are necessary, this method could be used as the first step of iterative methods to speed up convergence.} Last the method has been extended to any dimension $d$ and amounts to inverting an $N^d\times N^d$ Toeplitz matrix. The method \msgguy{can in particular be applied to images}.


\section*{Acknowledgments}
GP wishes to thank Sylvain Durand for his very important advices to write the paper and the anonymous referee who helped to improve the paper.


%
%
%

\bibliographystyle{model1-num-names}
\bibliography{references}

%

\end{document}